\documentclass{article} %\documentclass[11pt]{article}
\usepackage{amsthm,amsmath,amssymb,mathrsfs,amscd,color}
\newtheorem{thm}{Theorem}[section]
\newtheorem{cor}[thm]{Corollary}
\newtheorem{lem}[thm]{Lemma}
\newtheorem{defn}[thm]{Definition}

\newenvironment{main}[1]{\par\medskip\noindent{\bf #1.}\itshape}{\par\medskip}
\DeclareMathOperator{\ord}{ord}
\DeclareMathOperator{\Aut}{Aut}

\DeclareMathOperator{\Irr}{Irr}
\DeclareMathOperator{\IBr}{IBr}
\DeclareMathOperator{\Hom}{Hom}
\DeclareMathOperator{\Gal}{Gal}
\begin{document}
\title{Hyperbolic Modules of Finite Group Algebras 
  over Finite Fields of Characteristic Two}
\author{Ping Jin\thanks{Supported by the NSF of China (No. 11171194)
and partially by the NSF of China (No. 11071091).}\\
\small School of Mathematical Sciences, Shanxi University, Taiyuan 030006, PR China\\
\small{\it Email address}: jinping@sxu.edu.cn~\\ [5pt]
Yun Fan\thanks{Supported by the NSF of China (No. 11271005).}\\
\small School of Mathematics and Statistics, Central China Normal University,\\
\small Wuhan 430079, PR China\\
\small{\it Email address}: yfan@mail.ccnu.edu.cn}
\date{}
\maketitle

%%%----------------------------------------------------------------------
\begin{abstract}
Let $G$ be a finite group and let $F$ be a finite field of
characteristic~$2$. We introduce \emph{$F$-special subgroups} and
\emph{$F$-special elements} of $G$.
In the case where $F$ contains a $p$th primitive root of unity
for each odd prime $p$ dividing the order of $G$
(e.g. it is the case once $F$ is a splitting field for all subgroups of $G$),
the $F$-special elements of $G$ coincide with real elements of odd order.
We prove that a symmetric $F G$-module $V$ is hyperbolic if and only if the
restriction $V_D$ of $V$ to every $F$-special subgroup $D$ of $G$ is
hyperbolic, and also, if and only if the characteristic polynomial on
$V$ defined by every $F$-special element of $G$ is a square of a
polynomial over $F$.
Some immediate applications to characters, self-dual codes and Witt groups are given.


\end{abstract}

%%%%%%%%%%%%%%%%%%%%%%%%%%%%%%%%%%%%%%%%%%%%%%%%%%%%%%%%%%%%%%%%%%%%%%%%%%%%%%%%%%%%%%%%%%
\section{Introduction}
Let $G$ be a finite group and $F$ be a finite field. Let $V$ be a
finite-dimensional $FG$-module which carries
a $G$-invariant non-degenerate bilinear form.
We call $V$ a \emph{symplectic} (\emph{symmetric} respectively)
$FG$-module if the form is alternating (symmetric respectively).

Let $V$ be a symmetric or symplectic $FG$-module.
A submodule $W$ of $V$ is said to be {\em isotropic} if $W\subseteq W^\bot$;
while $W$ is said to be {\em self-perpendicular} if  $W^\bot=W$.
Following Dade~\cite{D}, we say that $V$ is \emph{hyperbolic}
if it contains a  self-perpendicular submodule $W$.
The hyperbolic modules are also named \emph{metabolic} modules in many literatures.

The theory on hyperbolic modules is powerful and fruitful in many areas.
The hyperbolic symplectic modules are quite useful for investigating
correspondences of characters, extensions of characters,
monomial characters and $M$-groups, for example,
see~\cite{D,I1983,Le,L,L2006,P,Wa,WH}.
The study of hyperbolic symmetric modules is used to determine
the existence of self-dual codes and to characterize
the automorphisms of self-dual codes,  see~\cite{FJ,FZ,GN,MW,MW3}.
Also, in order to study the Frobenius-Schur indicators of group representations,
J.G. Thompson~\cite{T} defined the \emph{Witt kernel} of a symmetric or symplectic $FG$-module $V$
to be $W^\bot/W$ for a maximal isotropic submodule $W$ of $V$,
which inherits the form on $V$ and is unique up to $G$-isometry
(see~\cite[Lemma 2.1]{T}).
Obviously, $V$ is hyperbolic if and only if its Witt kernel is trivial.

Dade \cite[Theorem 3.2]{D} proved the following theorem
which reduces the problem of determining whether or not
a symplectic module is hyperbolic to the ``semisimple case''.

\begin{thm}[Dade]
Suppose that $F$ is a finite field of odd characteristic $p$,
that $G$ is a $p$-solvable group, that $H$ is a subgroup of $p$-power
index in $G$, and that $V$ is a symplectic $F G$-module whose restriction $V_H$
to a symplectic $F H$-module is hyperbolic. Then $V$ is hyperbolic.
\end{thm}

Unfortunately, the theorem does not hold  for characteristic $p=2$,
even if $F$ is a splitting  field for all subgroups of $G$.
A counterexample can be constructed over $\mathbb{F}_4[S_3]$,
where $\mathbb{F}_4$ denotes the field with four elements and
$S_3$ denotes the symmetric group on three letters.
Not surprisingly, ``characteristic $2$" causes great difficulties in many problems.

Using the above Dade's theorem, Loukaki~\cite[Theorem A]{L}
present an effective criterion for finite groups of odd order.

\begin{thm}[Loukaki]
Suppose that $F$ is a finite field of odd characteristic $p$,
that $G$ is a finite group of odd order, and that $V$ is a symplectic $F G$-module
whose restriction $V_C$ to a symplectic $F C$-module is hyperbolic
for every cyclic subgroup $C$ of $G$. Then $V$ is hyperbolic.
\end{thm}

In the situation of the above Loukaki' theorem,
we removed the assumption that the characteristic is odd,
reduced the criterion to fewer elements,
and showed one more criterion by characteristic polynomials.
Moreover, our proof is no longer based on Dade's theorem,
see Corollaries B and C in~\cite{J2011}.

However, it seems more difficult and complicated to handle with the case where
both  the characteristic of the ground field and the order of the group are even.
Dade's arguments used in~\cite{D} and ours in~\cite{J2011}
are no longer valid in this ``modular" case.

The main purpose of the present paper is to explore general and effective criteria to
determine the hyperbolicity of symmetric or symplectic modules
in the case of characteristic two.
In order to obtain a similar criterion by characteristic polynomials,
a key idea is how to carefully pick out some distinguished elements
that can be used to determine if a symmetric or symplectic module is hyperbolic.

To state our main result, we introduce necessary notations.
Fix a finite field $F$ with $q$ elements where $q$ is a power of $2$.
For each $i\geq0$, we write $\pi_i(q)$ for the set of
the odd primes $p$ satisfying that
$\nu_2(\ord_p(q))=i$, where $\ord_p(q)$ denotes the order of $q$ modulo $p$,
and $\nu_2$ is the $2$-adic discrete valuation on the rational field
$\mathbb{Q}$; i.e.
an odd prime $p\in\pi_i(q)$ if and only if $\ord_p(q)=2^ia$ with $2\nmid a$.
As usual, by $C_m$ we denote the cyclic group of order $m$,
by $D_{2n}$ we denote the dihedral group of order $2n$ with $n\geq3$.
In addition, for any integers $e\geq1$ and $n\geq3$ with $2\nmid n$,
we define the \emph{extended dihedral group} of order $2^en$ as
$$\tilde D_{2^en}=\langle x,y\,|\,x^n=y^{2^e}=1,\,y^{-1}xy=x^{-1}\rangle.$$
Clearly $\tilde D_{2^en}/O_2(\tilde D_{2^en})\cong D_{2n}$,
and $\tilde D_{2^en}=D_{2n}$ for $e=1$.
Here $O_2(G)$ denotes the maximal normal $p$-subgroup of a finite group $G$.

We define the \emph{$F$-special groups} and \emph{$F$-special elements}
as follows, which play a key role in our study of hyperbolic modules.

\begin{defn}\rm
Let $F$ be a finite field of even order $q$.
A finite group $G$ is said to be {\em $F$-special},
if it has one of the following isomorphism types:

{\rm(i)} $C_m$, where $m\geq1$ is a $\pi_i(q)$-number for some $i\geq1$;

{\rm(ii)} $\tilde D_{2^en}$;

{\rm(iii)} $C_m\times \tilde D_{2^en}$, where $m\geq3$ is a $\pi_i(q)$-number for some $i\geq1$,
and $n$ is odd with $\nu_2(\ord_n(q))<i$.
Note that $\gcd(m,n)=1$ in this case and hence $G\cong C_{mn}\rtimes C_{2^e}$.

Similarly, an element $g$ of  a finite group $G$ is called {\em $F$-special},
if it has odd order and satisfies one of the following (possibly overlapped) conditions:

{\rm(1)} $g$ is a $\pi_i(q)$-element for some $i\geq1$;

{\rm(2)} $g$ is real, that is, $g^h:=h^{-1}gh=g^{-1}$ for some $h\in G$;

{\rm(3)} There exists an involution (elements of order $2$)
$$t\in N_G(\langle g\rangle)/C_G(g)$$
such that
$C_{\langle g\rangle}(t)$ is a nontrivial $\pi_i(q)$-group for some $i\geq1$,
and $\nu_2(\ord_n(q))<i$ where $n$ denotes the order of the commutator $[g,t]=g^{-1}g^t$.
\end{defn}

By the definition, it is clear that if the finite field $F$ contains
a $p$th primitive root of unity for each odd prime $p$ dividing the order $|G|$ of a finite group $G$
(e.g. it is the case once $F$ is a splitting field for all subgroups of $G$),
then every odd prime divisor of $|G|$ belongs to $\pi_0(q)$,
which implies that the $F$-special nontrivial subgroups of $G$ coincide with
the extended dihedral subgroups,
and that an element of $G$ is $F$-special if and only if it is real of odd order.
(Note that for any real element $g$ of $G$, one can choose a $2$-element $t\in G$
such that $g^t=g^{-1}$.)

Now we state our main results.
Clearly, we only need to deal with symmetric modules,
since symplectic modules over any field of characteristic two
are automatically symmetric.

\begin{main}{Theorem A}
Let $G$ be a finite group, $F$ be a finite field of characteristic~$2$
and $V$ be a symmetric $FG$-module.
Then the following three statements are equivalent.

{\rm(1)} $V$ is hyperbolic.

{\rm(2)} $V_D$ is hyperbolic for every $F$-special subgroup $D$ of $G$.

{\rm(3)} The characteristic polynomial on $V$ defined by every $F$-special element of $G$
is a square of a polynomial over $F$.

In particular, if $F$ contains a $p$th primitive root of unity
for each odd prime $p$ dividing $|G|$,
then $V$ is hyperbolic if and only if the characteristic polynomial on $V$
defined by every real element of $G$ of odd order is a square.
\end{main}

We remark that, though the definition of $F$-special groups
and $F$-special elements seems out of ordinary,
it plays an essential role in determining whether or not
a symmetric module is hyperbolic.
We will present examples in Section 5 to show that, in some sense,
the set of the $F$-special subgroups of a finite group is
``minimal" to guarantee that Theorem A holds.
Another remark is that, of course, \emph{not} every element of $G$
has a square characteristic polynomial in general,
even if $V$ is hyperbolic and $G$ is cyclic,
as shown by~\cite[Example 2.5]{J2011}.

To prove Theorem A, we will work with self-dual $FG$-modules.
Recall that a finite-dimensional
$FG$-module $V$ is said to be \emph{self-dual}
if $V$ is isomorphic to its dual $V^*=\Hom_F(V,F)$ as $FG$-modules.
An obvious fact is that an $FG$-module is self-dual if and only if
it carries a $G$-invariant non-degenerate bilinear form (see~\cite[VII.8.10]{H2}).
In particular, every symmetric $FG$-module is self-dual,
but a self-dual module in characteristic $2$
can fail to be a symmetric module, as shown by an example in~\cite[page 107]{H2}.
Based on a Dade's result~\cite[Corollary 2.10]{D},
it easily follows that a symmetric module in characteristic two
is hyperbolic if and only if every self-dual composition factor occurs
with even multiplicity (see Lemma~\ref{em} below).
We may therefore view Theorem A as an immediate corollary of the following Theorem B
that is what we really want to prove.

For convenience, we say that a self-dual $FG$-module $V$
is an \emph{even-multiplicity} module if every
self-dual simple $FG$-module occurs with even multiplicity
in a composition series of $V$.

\begin{main}{Theorem B}
Let $G$ be a finite group, $F$ be a finite field of characteristic $2$
and $V$ be a self-dual $FG$-module.
Then the following three statements are equivalent.

{\rm(1)} $V$ is an even-multiplicity $FG$-module.

{\rm(2)} 
$V_D$ is an even-multiplicity $FD$-module for every $F$-special subgroup $D$ of~$G$.

{\rm(3)} The characteristic polynomial on $V$ defined by every $F$-special element of $G$
is a square of a polynomial over $F$.

In particular, if $F$ contains a $p$th primitive root of unity
for each odd prime $p$ dividing $|G|$,
then $V$ is an even-multiplicity module if and only if the characteristic polynomial on $V$
defined by every real element of $G$ of odd order is a square.
\end{main}

The proof of Theorem B constitutes a big part of this article.
The most crucial step towards the proof of Theorem B is the following
Theorem C, which looks concise and we think is of
independent interest. Following~\cite[Definition 8.8]{I}, a finite
group is said to be \emph{$2$-quasi-elementary} if it has the form
$C\rtimes P$, where $C$ is a cyclic group of odd order acted on by a
$2$-group $P$. Note that the ground field $F$ involved in Theorem C
can be arbitrary, not necessarily finite.

\begin{main}{Theorem C}
Let $G$ be a finite group and $F$ be any field.
Suppose that $V$ is a self-dual $FG$-module
whose restriction $V_H$ to every $2$-quasi-elementary subgroup $H$ of $G$
is an even-multiplicity $FH$-module.
Then $V$ is an even-multiplicity $FG$-module.
\end{main}

We remark that Loukaki's theorem mentioned above
is an immediate consequence of Theorem C
(along with an elementary fact that any symplectic module of finite groups of odd order
is hyperbolic if and only if it is an even-multiplicity module,
see Corollary 3.8 in~\cite{J2011}).
We may also use this theorem to study the real-valued complex characters
or $2$-Brauer characters of a finite group
by taking $F$ to be a splitting field for $G$ of characteristic zero or two
(see Theorem~\ref{real-valued Brauer} below).

The organization of the paper is as follows.
In Section 2 we prove Theorem C by virtue of
the Grothendieck group $G_0(FG)$ of the category of all finitely generated $FG$-modules
and the version of the Solomon induction theorem~\cite[Theorem 8.11]{I} over arbitrary field.
Clearly, Theorem C allows us to prove Theorem B
only for $2$-quasi-elementary groups.
In Section 3 we study in detail the self-duality of simple modules of
$2$-quasi-elementary groups of special types,
and then we prove that the $F$-special groups defined above are
just the groups whose simple modules are all self-dual.
Section 4 is devoted to proving Theorem B and Theorem A.
Section 5 provides two examples of $F$-special subgroups.
In Section 6 we give some immediate applications of our theorems to
the real-valued irreducible Brauer characters at prime $2$,
the automorphisms of self-dual binary codes,
and the Witt group of symmetric or symplectic modules.

Throughout this paper, all groups are finite
and all modules are finite-dimensional right modules.
Notations and conventions are mainly from Huppert's book~\cite{H2} and from Isaacs' book~\cite{I}.
Also, we will frequently use the Clifford theorem over arbitrary fields (see~\cite[Theorem 6.5]{I})
and the Krull-Schmidt theorem for modules (see~\cite[I.12.4]{Hun}),
without further explicit references.

%%%%%%%%%%%%%%%%%%%%%%%%%%%%%%%%%%%%%%%%%%%%%%%%%%%%%%%%%%%%%%%%%%%%%%%%%%%%%%%%%%%%%%%%
\section{Proof of Theorem C}
In this section we will prove Theorem C.
To do this, we need to work with the Grothendieck group $G_0(FG)$
of the category of all finitely generated $FG$-modules.
For definition and properties, see~\cite{CR,Se}, for example.

By a \emph{$p$-modular system} for a fixed prime $p$, we mean a triple $(K,R,F)$
where $R$ is a complete discrete valuation ring of characteristic $0$,
$K$ is the quotient field of $R$, and $F$ is the residue field of $R$ of characteristic $p$.

We begin with a useful characterization of even-multiplicity modules
over any field in terms of Grothendieck groups.

\begin{lem}
Let $G$ be any group and let $F$ be any field.
Suppose that $V$ is a self-dual $FG$-module.
Then $V$ is an even-multiplicity module if and only if
$$[V]=2[\Delta]+[\Lambda]+[\Lambda^*]$$
in $G_0(FG)$ for some $FG$-modules $\Delta$ and $\Lambda$,
where $\Lambda^*$ denotes the dual of $\Lambda$.
\end{lem}
\begin{proof}
Clear from definition.
\end{proof}

Next is a version of the Solomon's induction theorem
(see~\cite[Theorem 8.10]{I} or ~\cite{So}) over arbitrary fields.
For any subgroup $H$ of $G$, we denote by $F(G/H)$
the permutation $FG$-module defined by the $G$-set $G/H=\{Hg\,|\,g\in G\}$
that is the set of right cosets of $H$ in $G$ on which $G$ acts by multiplication.

\begin{lem}
Let $F$ be any field and let $G$ be any group.
Then there exist some $2$-quasi-elementary subgroups $H_1,\ldots,H_n$ of $G$,
some integers $a_1,\ldots,a_n$, and an odd integer $m$,
such that for any $FG$-module $V$, we have
$$ m[V]=\sum_{i=1}^na_i[(V_{H_i})^G]$$
in $G_0(FG)$.
\end{lem}
\begin{proof}
The Solomon's induction theorem mentioned above tells us that
$$m[\mathbb{Q}]=a_1[\mathbb{Q}(G/H_1)]+\cdots+ a_n[\mathbb{Q}(G/H_n)] \eqno{(2.1)}$$
in $G_0(\mathbb{Q}G)$, where each $H_i$, $1\leq i\leq n$,
is a $2$-quasi-elementary subgroup of $G$, $a_i$ is an integer,
$m$ is odd, and the rational field $\mathbb{Q}$ is the trivial $\mathbb{Q}G$-module.
We will prove that Equation (2.1) remains true for any field $F$
in the place of $\mathbb{Q}$, that is,
$$m[F]=a_1[F(G/H_1)]+\cdots+ a_n[F(G/H_n)] \eqno{(2.2)}$$
in $G_0(FG)$.

Assume first that $F$ is of characteristic $0$.
Then we may let $\mathbb{Q}\subseteq F$ without any loss.
In this case, by scalar extension from $\mathbb{Q}$ to $F$,
we have an additive map $G_0(\mathbb{Q}G)\to G_0(FG)$,
which yields Equation (2.2).

Assume that $F$ has prime characteristic $p$.
Then by the MacLane's existence theorem (see Theorem 3.1.5 of~\cite{Kar}),
we may choose a $p$-modular system $(K,R,F)$ containing $F$.
Since $K$ is of characteristic $0$,
it follows that Equation (2.2) holds for $K$ in the place of $F$
by the same reason.
Moreover, it is easy to see that
$K(G/H_i)\cong R(G/H_i)\otimes_RK$ as $KG$-modules and
$F(G/H_i)\cong R(G/H_i)\otimes_RF$ as $FG$-modules for each $i=1,\ldots,n$,
which implies that the decomposition map
$$d:G_0(KG)\to G_0(FG)$$
sends $[K(G/H_i)]$ into $[F(G/H_i)]$ by definition.
In passage from $K$ to $F$ via
the map $d$, we conclude that Equation (2.2) also holds
for any field $F$ of prime characteristic.

Now we denote by $1_H$ the trivial $FH$-module $F$ for any subgroup $H$ of $G$.
Then we have $(1_H)^G\cong F(G/H)$ as $FG$-modules.
Let $V$ be an $FG$-module. By~\cite[Lemma VII.4.15]{H2}, we see that
$$V\otimes_F F(G/H_i)=V\otimes_F(1_{H_i})^G\cong (V_{H_i}\otimes_F1_{H_i})^G\cong(V_{H_i})^G$$
as $FG$-modules,
and the result follows by multiplying with $[V]\in G_0(FG)$ in both sides of Equation (2.2).
\end{proof}

The following is Theorem C, which we restate here for convenience.

\begin{thm}\label{2-quasi}
Let $G$ be any group and let $F$ be any field.
Suppose that $V$ is a self-dual $FG$-module
whose restriction $V_H$ to every $2$-quasi-elementary subgroup $H$ of $G$
is an even-multiplicity $FH$-module.
Then $V$ is an even-multiplicity $FG$-module.
\end{thm}
\begin{proof}
By Lemma 3.2, we see that in $G_0(FG)$,
$$ m[V]=a_1[(V_{H_1})^G]+\cdots+a_n[(V_{H_n})^G],$$
where all $H_i$ are $2$-quasi-elementary subgroups of $G$,
$a_i$ are integers and $m$ is odd.
Since each restriction $V_{H_i}$ is an even-multiplicity $FH_i$-module,
it follows from Lemma 3.1 that
$[V_{H_i}]=2[\Delta_i]+[\Lambda_i]+[\Lambda_i^*]$
in $G_0(FH_i)$ for some $FH_i$-modules $\Delta_i$ and $\Lambda_i$,
and hence we have
$$m[V]=\sum_{i=1}^n a_i\left(2[\Delta_i^G]+[\Lambda_i^G]+[(\Lambda_i^*)^G]\right).\eqno{(2.3)}$$
For any self-dual simple $FG$-module $U$,
it is easy to see that $U$ has the same multiplicity
as a composition factor in both $\Lambda_i^G$
and its dual $(\Lambda_i^G)^*\cong (\Lambda_i^*)^G$.
Hence $U$ must occur with even multiplicity as a composition factor of $V$ by Equation (2.3)
and $V$ is an even-multiplicity $FG$-module, as desired.
\end{proof}

%%%%%%%%%%%%%%%%%%%%%%%%%%%%%%%%%%%%%%%%%%%%%%%%%%%%%%%%%%%%%%%%%%%%%%%%%%%%%%%%%%%%%%%%%%%%%%%%%%%%%%
\section{$2$-quasi-elementary groups}
Throughout this section $F$ always denotes a finite field of characteristic $2$ with $q$ elements,
unless otherwise stated.
One of our aims is to prove that every simple module of $F$-special groups
is self-dual.

We begin with Green's theorem that will be used frequently in the present paper
without further reference.

\begin{lem}
Let $F$ be a finite field of prime characteristic $p$,
and let $N$ be a normal subgroup of a group $G$ with $p$-power index.
If $W$ is a simple $FN$-module,
then there exists a unique, up to isomorphism, simple $FG$-module $V$ lying over $W$.
Furthermore, if $W$ is $G$-invariant, then $V$ is an extension of $W$ to $G$.
\end{lem}
\begin{proof}
See~\cite[Lemma VII.9.19]{H2}, for example.
\end{proof}

As a preliminary we need an elementary result on faithful simple
modules for $2$-quasi-elementary groups.

\begin{lem}\label{faith-mod}
Let $G$ be a $2$-quasi-elementary group with $O_2(G)=1$,
and let $N$ be the normal $2$-complement of $G$.
Suppose that $V$ is a simple $FG$-module and that $W$ is a simple $FN$-submodule of $V$.
Then $V$ is a faithful $FG$-module if and only if $W$ is a faithful $FN$-module.
\end{lem}
\begin{proof} Since $N$ is cyclic and normal in $G$,
we see that the kernel $C_N(W)$ of the action of $N$ on $W$ is also normal in $G$.
It easily follows from Clifford's theorem that
$$C_G(V)\cap N=C_N(V)=\bigcap_{g\in G} C_N(W^g)=C_N(W).$$

If $V$ is a faithful $FG$-module, that is, $C_G(V)=1$,
then $C_N(W)=1$ and hence $W$ is also faithful as a $FN$-module.

Conversely, if $W$ is faithful, then $C_G(V)\cap N=1$,
which implies that $C_G(V)$ is a normal $2$-subgroup of $G$ and contained in $O_2(G)$.
So $C_G(V)=1$ and $V$ is a faithful $FG$-module.
\end{proof}

Now, for simplicity we introduce an integer-valued function $\omega$,
defined on odd integers, groups and modules, respectively,
with respect to the fixed finite field $F$ of characteristic $2$ with $q$ elements.

$\bullet$ For any odd integer $n\geq1$, define $\omega(n)=\nu_2(\ord_n(q))$.
It is easy to prove that
$$\omega(n)=\max\{\omega(p)\;|\;p\; \text{is a prime divisor of}\; n\}.$$
Moreover, for any odd prime $p$ and any integer $i\geq0$,
by definition we see that
$$p\in\pi_i(q)\Longleftrightarrow \omega(p)=i.$$

$\bullet$ For any finite group $G$, define $\omega(G)=\max\{\omega(p)\}$
where $p$ runs over all odd prime divisors of $|G|$.
We also have $\omega(G)=\omega(|G|_{2'})$.

$\bullet$ For any finite dimensional $F$-space $V$, define $\omega(V)=\nu_2(\dim V)$.

\medskip
A simple property on the function $\omega$ is as follows,
which will be used in the proof of Theorem B.
\begin{lem}\label{omega}
Let $G$ be a group, and let $V$ be a faithful simple $F G$-module.

{\rm (1)} If $G$ is cyclic of odd order, then $\omega(V)=\omega(G)$.

{\rm (2)} If $G$ is dihedral with $4\nmid |G|$, then $\omega(V)\leq\omega(G)+1$.
\end{lem}
\begin{proof}
(1) Let $|G|=n$ with $2\nmid n$.
Then $\dim V=\ord_n(q)$ by~\cite[Problem 9.20]{I},
and hence $\omega(V)=\omega(G)$ by definition.

(2) Let $G=N\rtimes\langle t\rangle$, where $N$ is cyclic of odd order $n$, and
$t$ is an involution and inverts $N$.
By Clifford's theorem, we see that $V_N\cong W$ or $W\oplus W^t$
for some simple $F N$-module $W$.
Hence $\dim V=\dim W$ or $\dim V=2\dim W$, which implies that $\omega(V)\leq\omega(W)+1$.
Since $V$ is a faithful $FG$-module,
it follows that $O_2(G)=1$.
Applying Lemma~\ref{faith-mod}, we know that
$W$ is also a faithful $FN$-module.
So $\omega(W)=\omega(N)=\omega(G)$ by (1),
which yields that $\omega(V)\leq\omega(G)+1$.
\end{proof}

We now examine the duality of simple modules for $F$-special groups.
The ``cyclic case" is quite simple and already appeared in~\cite{J2011}.
\begin{lem}\label{cyc-dual}
Let $G$ be a cyclic group of odd order. Then the following assertions are equivalent.

{\rm(1)} Each simple $F G$-module is self-dual.

{\rm(2)} There exists a faithful self-dual simple $F G$-module.

{\rm(3)} $G$ is a $\pi_i(q)$-group for some $i\geq1$.
\end{lem}
\begin{proof}
See Lemmas 2.1 and 2.2 in~\cite{J2011}.
\end{proof}

\begin{lem}\label{key}
Let $G=N\rtimes\langle t\rangle$ be a non-abelian group
where $N$ is cyclic of odd order and $t$ is an involution,
and suppose that $W$ is a faithful simple $F N$-module.
Then $W^t\cong W^*$ as $F N$-modules
if and only if one of the following conditions holds:

{\rm(1)} $G\cong D_{2n}$, where $n$ is odd.

{\rm(2)} $G\cong C_m\times D_{2n}$, where $m\geq3$ is a $\pi_i(q)$-number for some $i\geq1$,
and $n$ is odd with $\omega(n)<i$.
\end{lem}
\begin{proof}
Let $N=\langle g\rangle$ be of odd order $r$, and let $g^t=g^k$ for some integer $k$.
By hypothesis we conclude that $k\not\equiv1\;({\rm mod}\; r)$ and $k^2\equiv1\;({\rm mod}\; r)$.
Let $c(x)$ be the characteristic polynomial of $g$ on $W$.
Then, by Lemma 2.1 of~\cite{J2011}, we see that $c(x)$ is irreducible over $F$
and thus we may assume that $\xi,\xi^q,\ldots,\xi^{q^{d-1}}$
are all roots of $c(x)$ in its a splitting field,
where $\xi$ is a $r$-th primitive root of unity and $d$ is the $F$-dimension of $W$
which also equals to $\ord_r(q)$.
Furthermore, we may easily deduce that $W^t\cong W^*$ if and only if
$g^t$ and $g^{-1}$ have the same characteristic polynomial on $W$,
and also if and only if $\xi^{kq^i}=\xi^{-1}$ for some positive integer $i$,
that is, the congruence equation
$$kq^x\equiv-1\;({\rm mod}\; r)\eqno{(3.1)}$$
has a positive integer solution $x$.
Note that if $k\equiv-1\;({\rm mod}\; r)$,
that is, if $G$ is dihedral,
then the above equation always has a solution,
and the result follows in this case.
So we may assume further that $k\not\equiv-1\;({\rm mod}\; r)$.
In this situation we may write $r=mn$ with $m,n\geq3$, $m|(k-1)$ and $n|(k+1)$.
Then $t$ centralizes $g^n$ and inverts $g^m$, which yields that
$$G=\langle g\rangle\rtimes\langle t\rangle=
(\langle g^n\rangle\times\langle g^m\rangle)\rtimes\langle t\rangle
=\langle g^n\rangle\times(\langle g^m\rangle\rtimes\langle t\rangle)
\cong C_m\times D_{2n}.$$
What remains is to prove that $W^t\cong W^*$ if and only if condition (2) holds.

Suppose first that $W^t\cong W^*$, that is, Equation (3.1) has a positive integer solution $x$.
Then $q^x\equiv-1\;({\rm mod}\; m)$. Let $\nu_2(x)=i-1$ for some $i\geq1$.
It is easy to deduce that $m$ is a $\pi_i(q)$-number.
Also, we see that $q^x\equiv1\;({\rm mod}\; p)$ for every prime $p$ dividing $n$,
and thus the order $\ord_p(q)$ of $q$ modulo $p$ must divide $x$.
It follows that $\omega(p)\leq i-1$ and hence $\omega(n)<i$.
This proves that $G$ satisfies condition (2), as wanted.

Conversely, suppose that condition (2) holds and we only need to prove that
Equation (3.1) has a positive solution.
In this case, the hypothesis on $m$ easily implies that there exists an odd integer $s$
satisfying $q^{2^{i-1}s}\equiv-1\;({\rm mod}\; m)$.
(See the last paragraph in the proof of~\cite[Lemma 2.2]{J2011}, for example.)
Also, since $\omega(n)<i$, it follows that the $2$-part of $\ord_n(q)$
is less than $2^i$ and hence we have
$$q^{{2^{i-1}}s'}\equiv1\;({\rm mod}\; n)$$
for some odd integer $s'$. Therefore Equation (3.1) has a solution $x=2^{i-1}ss'$,
which implies that $W^t\cong W^*$. The proof is now complete.
\end{proof}

For convenience, we introduce a subclass of $2$-quasi-elementary groups,
which contains all cyclic groups of odd order and all $F$-special groups.

\begin{defn}\rm
A group $G$ is said to be \emph{hyper-dihedral}
if $G=\langle x\rangle\rtimes\langle t\rangle$,
where $x$ is of odd order and $t$ is a $2$-element,
such that $t^2$ centralizes $x$ but $t$ does not, unless $t=1$.
\end{defn}

We will determine the duality of simple modules for such groups.
Since each simple $FG$-module can be viewed as a simple $F(G/O_2(G))$-module,
we may assume that $O_2(G)=1$ without any lose.
By Lemma~\ref{cyc-dual}, we need only consider non-abelian hyper-dihedral groups.

\begin{thm}\label{dih-dual}
Let $G$ be a non-abelian hyper-dihedral group with  $O_2(G)=1$.
Then the following assertions are equivalent.

{\rm(1)} Each simple $FG$-module is self-dual.

{\rm(2)} There exists a faithful self-dual simple $FG$-module.

{\rm(3)} $G$ is isomorphic to one of the following groups:

\ \quad {\rm(i)} $D_{2n}$, where $n$ is odd;

\ \quad {\rm(ii)} $C_m\times D_{2n}$, where $m\geq3$ is a $\pi_i(q)$-number for some $i\geq1$,
and $n$ is odd satisfying $\omega(n)<i$;

\ \quad {\rm(iii)} $C_m\times D_{2n}$, where $m\geq3,\,\gcd(m,n)=1$,
and $mn$ is a $\pi_i(q)$-number for some $i\geq1$.
\end{thm}
\begin{proof}
By hypothesis, we may let $G=N\rtimes\langle t\rangle$,
where $N$ is cyclic of odd order and $t$ is an involution.
Then $N=C_N(t)\times[N,t]$ by Fitting's lemma.
Also, since $G$ is non-abelian, it follows that $G\cong C_m\times D_{2n}$
in which both $m$ and $n$ are odd.

That ``(1)$\Rightarrow$(2)" is clear by Lemma~\ref{faith-mod},
since the cyclic group $N$ always has a faithful simple module.

 ``(2)$\Rightarrow$(3)".
If $m=1$, then $G\cong D_{2n}$ and result follows.
So we may assume that $m>1$.
Let $V$ be a faithful self-dual simple $F G$-module,
and take $W$ to be a simple $F N$-submodule of $V$.
Then $W$ is also faithful by Lemma~\ref{faith-mod}.
If $W^t\cong W^*$ then $G$ is of isomorphism type (i) or (ii) by Lemma~\ref{key}.
And if $W^t\ncong W^*$ then $W$ must be self-dual by Clifford's theorem.
In this case, $N$ is a $\pi_i(q)$-group for some $i\geq1$ by Lemma~\ref{cyc-dual},
which implies that $G$ has isomorphism type (iii), as wanted.

 ``(3)$\Rightarrow$(1)".
Let $V$ be any simple $FG$-module and set $K=C_G(V)$.
Then $V$ is a faithful simple $G/K$-module over $F$ so that $O_2(G/K)=1$.
If $t\in K$, then $[N,t]\subseteq K$,
which implies that $G/K$ must be a cyclic $\pi_i(q)$-group for some $i\geq1$
by the isomorphic types of $G$ listed in (3).
Hence $V$ is self-dual by Lemma~\ref{cyc-dual} and we are done in this case.
So we may assume further that $t\not\in K$.
Then we have $K\subseteq N$. Clearly $K\neq N$ as $O_2(G/K)=1$ and
thus $G/K\cong N/K\rtimes\langle t\rangle$ is also a non-abelian hyper-dihedral group
whose isomorphism type is presented by (3).
For simplicity we may let $K=1$ so that $V$ is faithful.
By Lemma~\ref{faith-mod} again, we see that $V$ has a faithful simple $FN$-submodule $W$.
If $G$ is of type (i) or (ii), then Lemma~\ref{key}
implies that $W^t\cong W^*$.
And if $G$ is of type (iii) then $W\cong W^*$ by Lemma~\ref{cyc-dual}.
Hence we have $V\cong V^*$ in all cases by Green's theorem,
which proves (1). The proof is now complete.
\end{proof}

Combining Lemma~\ref{cyc-dual} with Theorem~\ref{dih-dual}, we have
\begin{cor}\label{special-dual}
If $G$ is an $F$-special group, then every simple $FG$-module is self-dual.
\end{cor}

We end this section with a remark.
By Lemma~\ref{cyc-dual} and Theorem~\ref{dih-dual}, we see that $F$-special groups $D$
are essentially those hyper-dihedral groups all of whose simple modules are self-dual,
with an exception that $D/O_2(D)$ is of type (iii) in Theorem~\ref{dih-dual}.
However, this exceptional case is superfluous for the proof of Theorem B,
since it can be included in cyclic $F$-special subgroups by Green's theorem.

%%%%%%%%%%%%%%%%%%%%%%%%%%%%%%%%%%%%%%%%%%%%%%%%%%%%%%%%%%%%%%%%%%%%%%%%%%%%%%%%%%%%%%%%
\section{Proofs of Theorems A and B}
In this section we will prove Theorems A and B from the Introduction.
We first establish the implication ``(2)$\Longrightarrow(1)$"
in Theorem B, which is the hard case of that theorem.

\begin{thm}\label{1}
Suppose that $G$ is a group,
that $F$ is a finite field of characteristic $2$ with $q$ elements,
and that $V$ is a self-dual $F G$-module
whose restriction $V_D$ to every $F$-special subgroup $D$ of $G$
is an even-multiplicity $F D$-module.
Then $V$ is an even-multiplicity $F G$-module.
\end{thm}

The proof of this theorem will follow from a series of lemmas,
all based on the hypothesis that $F,G,V$ form a minimal counterexample,
that is, we assume the following

\begin{main}{Inductive Hypothesis}
$F,G,V$ have been chosen among all the triplets satisfying the hypothesis
but not the conclusion of Theorem~\ref{1} so as to minimize:
first the order $|G|$ of $G$ and then the $F$-dimension $\dim V$ of $V$.
In particular, $G$ cannot be $F$-special.
\end{main}

We will derive a contradiction from the Inductive Hypothesis
by first reducing the group $G$ to a $2$-quasi-elementary group,
then to a hyper-dihedral group, and finally to an $F$-special group,
thus proving Theorem~\ref{1}.

The first step is an immediate application of Theorem C.
\begin{lem}\label{2}
$G$ is a $2$-quasi-elementary group of even order,
so that $G=N\rtimes T$, where $N$ is a cyclic normal subgroup of odd order,
and $T\neq1$ is a Sylow $2$-subgroup of $G$.
\end{lem}
\begin{proof}
The minimality of $|G|$ clearly implies that the restriction
$V_H$ of $V$ to each proper subgroup $H$ of $G$ is an even-multiplicity $FH$-module.
This, along with Theorem C from the Introduction, forces $G$ to be a $2$-quasi-elementary group.
So we may write $G=N\rtimes T$ for an odd-order cyclic group $N$
acted on by some $2$-group $T$.
If $T=1$, then $G$ is cyclic of odd order, and by~\cite[Theorem 2.4]{J2011},
we see that $G$ must be a cyclic $\pi_i(q)$-group for some $i\geq1$,
which cannot be the case. The proof is now complete.
\end{proof}

Next we will further reduce $G$ to a hyper-dihedral group, by Lemmas 4.3-4.9.
\begin{lem}\label{3}
We may assume that $$V\cong V_1\oplus\cdots\oplus V_n,\;n\geq1,$$
where $V_1,\ldots,V_n$ are pairwise non-isomorphic self-dual simple $F G$-modules.
\end{lem}
\begin{proof}
Let $V_1,\ldots,V_n$ be all the self-dual composition factors of $V$
that occur with odd multiplicity and let $V'=V_1\oplus\cdots\oplus V_n$.
It easily follows that, for any subgroup $H$ of $G$,
$V_H$ is an even-multiplicity $F H$-module if and only if $V'_H$ is.
So without any loss we may replace $V'$ by $V$ to find a desired contradiction.
\end{proof}

The following result shows that the above $V_i$'s are ``dually quasi-primitive".
\begin{lem}\label{4}
Let $L$ be a normal subgroup of $G$,
and let $X_i$ be an $FL$-submodule of $V_i$, for each $i=1,\ldots,n$.
Then each $FL$-submodule of $V_i$ is isomorphic with $X_i$ or its dual $X_i^*$.
In particular, $(V_i)_L\cong e_iX_i$ if $X_i$ is self-dual,
and $(V_i)_L\cong e_i(X_i\oplus X_i^*)$ otherwise, for some integer $e_i\geq1$.
\end{lem}
\begin{proof}
Fix some $i\in\{1,\ldots,n\}$. Set $U=V_i$ and $X=X_i$ for brevity.
Then Clifford's theorem implies that the restriction $U_L$
is a semisimple $FL$-module,
which contains $X$ and its dual $X^*$ as $FL$-submodules since $U$ is self-dual
by Lemma~\ref{3}. Let
$$W=U(X)+U(X^*)$$ be the sum of $X$-homogeneous part $U(X)$
and $X^*$-homogeneous part $U(X^*)$ of $U$.
Then we have $U_L=\sum_{g\in G}Wg$.
Write $J$ for the stabilizer of $W$ in $G$, that is, $J=\{g\in G\,|\,Wg=W\}$.
Then it is easy to see that the following three assertions hold:

(i) $W$ is a self-dual simple $FJ$-module.

(ii) $W^G\cong U$ as $FG$-modules.

(iii) $X$ has the same multiplicity in both $U_L$ and $W_L$,
which implies that $W$ has multiplicity $1$ as a composition factor of $U_J$.

Now suppose that $J$ is a proper subgroup of $G$.
Then $V_J$ is an even-multiplicity $FJ$-module by the minimality of $|G|$.
In this case, it follows from (i) and (iii) that
$W$ must occur as a composition factor in $(V_j)_J$ for some $j\neq i$.
Hence $X$ also lies under the $V_j$
and the same argument implies that $W^G\cong V_j$ by (ii),
which forces $V_i=U\cong V_j$, contradicting Lemma~\ref{3}.
So $J=G$ and $U_L=U(X)+U(X^*)$. The result follows.
\end{proof}

\begin{lem}\label{5}
Fix a maximal subgroup $M$ of $G$ containing $N$.
Then $M$ is normal in $G$ with $|G:M|=2$, and we have
$$(V_i)_M\cong U_i\oplus U_i^*,\quad i=1,\ldots n,$$
where $U_1,\ldots,U_n$ are pairwise non-isomorphic non-self-dual simple $FM$-modules.
\end{lem}
\begin{proof}
By Lemma~\ref{2} we see that $G/N\cong T$ is a nontrivial $2$-group,
which implies that $M/N$ is also maximal in $G/N$.
Hence $M$ is normal in $G$ with $|G:M|=2$.

Let $U_i$ be an $FM$-submodule of $V_i$ for each $i=1,\ldots,n$.
Since $V_M$ is an even-multiplicity $FM$-module
and each $V_i$ is self-dual,
it follows from Green's theorem that
all the $U_i$ are pairwise non-isomorphic simple $FM$-modules.
Furthermore each $U_i$ cannot be self-dual, because it has multiplicity $1$ in
the even-multiplicity module $V_M$.
Now Clifford's theorem yields that $(V_i)_M\cong U_i\oplus U_i^*$, as desired.
\end{proof}

\begin{lem}\label{6}
Set $W_i=(U_i)_N$ for each $i=1,\ldots,n$.
Then all $W_i$ are pairwise non-isomorphic non-self-dual simple $FN$-modules.
In particular, we have
$$(V_i)_N\cong W_i\oplus W_i^*$$
and $M=I_G(W_i)$, the inertia group of $W_i$ in $G$.
\end{lem}
\begin{proof}
This follows by Lemmas~\ref{4} and~\ref{5}, along with Green's theorem.
\end{proof}

\begin{lem}\label{7}
The Sylow $2$-subgroup $T$ of $G$ is cyclic.
\end{lem}
\begin{proof}
Since $T$ is nontrivial by Lemma~\ref{2}, we may let $S$ be any maximal subgroup of $T$.
Then $|T:S|=2$ and $NS$ is also a maximal subgroup of $G$ and contains $N$.
Now the same arguments used in Lemmas~\ref{5} and~\ref{6} will deduce that
$NS=I_G(W_i)=M$ for each $1\leq i\leq n$, which forces $S\subseteq M$
and $S=M\cap T$. Therefore $S$ is a unique maximal subgroup of $T$,
and by definition we know that $S$ must be the Frattini subgroup $\Phi(T)$ of $T$.
Now the fact that $|T:\Phi(T)|=2$ implies that $T$ is cyclic.
\end{proof}

We now investigate the faithfulness of the $FG$-module $V$.
\begin{lem}\label{faith}
$V$ is a faithful $FG$-module, that is, $C_G(V)=1$.
In particular, $O_2(G)=1$ and $G$ cannot be abelian.
\end{lem}
\begin{proof}
Let $K=C_G(V)$ be the kernel of the action of $G$ on $V$ and let $\bar G=G/K$.
Then $V$ can be viewed as an $F\bar G$-module.
It suffices to prove that $F,\bar G,V$ also satisfy the Inductive Hypothesis.
In this case, the minimality of $|G|$ will imply that $K=1$
and hence $V$ is a faithful $FG$-module.
Furthermore, since $V$ is semisimple by Lemma~\ref{3} and $F$ is of characteristic $2$,
it follows that $O_2(G)=1$, which implies that $G$ cannot be abelian as $T\neq1$ by Lemma~\ref{2}.

We fix an $F$-special subgroup $\bar D$ in $\bar G$,
and we only need to prove that $V_{\bar D}$ is an even-multiplicity $F\bar D$-module,
as $V$ cannot be an even-multiplicity $F\bar G$-module.
Note that if we can find some proper subgroup $H$ of $G$
such that the image $\bar H$ of $H$ in $\bar G$ coincides with $\bar D$,
then $V_{\bar D}$ must be an even-multiplicity module as $V_H$ is,
and we are done in this case.
So we may assume that no such proper subgroups $H$ exist.
From this we may deduce that $\bar D=\bar G$ and $K\subseteq\Phi(G)$
(the Frattini subgroup of $G$).
In what follows we will prove that $G$ itself is $F$-special,
which cannot be the case and thus completes the proof of the lemma.

Indeed, since both $N$ and $T$ are cyclic by Lemmas~\ref{2} and~\ref{7},
we may write $N=\langle g\rangle$ and $T=\langle t\rangle$.
Then $\bar D=\bar G=\langle\bar g\rangle\rtimes\langle\bar t\rangle$.
Also, since $V$ is a faithful semisimple $F\bar G$-module,
it follows that $O_2(\bar D)=O_2(\bar G)=1$ and hence $\bar t^2=1$
by definition of $F$-special groups. So $t^2\in K$.
Now the Fitting's lemma tells us that $N=[N,t^2]\times C_N(t^2)$
with $[N,t^2]\subseteq K\subseteq \Phi(G)$.
Hence $G=C_N(t^2)\rtimes\langle t\rangle$ and $N=C_N(t^2)$.
So $t^2$ centralizes $N$.
Of course, since $\bar G=\bar D$ is an $F$-special group with $O_2(\bar G)=1$,
we conclude that it must have one of the following simpler isomorphism types:

(a) $\bar G\cong C_m$, where $m$ is a $\pi_i(q)$-number for some $i\geq1$.

(b) $\bar G\cong D_{2n}$, where $n$ is odd.

(c) $\bar G\cong C_m\times D_{2n}$, where $m\geq3$ is a $\pi_i(q)$-number for some $i\geq1$,
and $n$ is odd with $\omega(n)<i$.

Assume (a). Then $\bar G=\langle\bar g\rangle$ and
hence $G=\langle g\rangle K=\langle g\rangle$,
which cannot be the case because $|G|$ is even by Lemma~\ref{2}.

Assume (b) or (c).
In both cases we see that $\bar N\cong C_m\times C_n$ (take $m=1$ in case (b)).
Since $N$ is cyclic and $\gcd(m,n)=1$, we may let $N=A\times B$,
where $A=\langle a\rangle$ and $B=\langle b\rangle$
such that $|\bar A|=m$ and $|\bar B|=n$.
Furthermore, by using the fact that $K\subseteq\Phi(G)$ again,
the elements $a,b\in N$ can be chosen so that
$|A|$ and $m$ have the same prime divisors (not counting multiplicity)
and that $|B|$ and $n$ also have the same prime divisors.
In this situation we see that $A$ is a cyclic $\pi_i(q)$-group as $m$ is a $\pi_i(q)$-number by (c).
Clearly $\bar t$ centralizes $\bar A$ and $A=[A,t]\times C_A(t)$
by the Fitting's lemma, so $[A,t]\subseteq K\subseteq\Phi(G)$.
From this we may deduce that $A=C_A(t)$ and hence $t$ centralizes $A$.
Also, since $B$ is of odd order, we may choose an element $c\in B$ such that $b=(c^2)^{-1}$.
Set $d=c^{-1}c^t$.
Then $\bar d=\bar c^{-2}=\bar b$ because $\bar t$ inverts $\bar B$.
Since we have proved that $t^2$ centralizes $N$,
it follows that $d^t=(c^t)^{{-1}}c^{t^2}=d^{-1}$
and thus $t$ normalizes the cyclic subgroup $\langle d\rangle$ of $B$.
This implies that $\langle a\rangle\times(\langle d\rangle\rtimes\langle t\rangle)$
is a subgroup of $G$ which maps onto $\bar G$ as
$\langle\bar d\rangle=\langle\bar b\rangle=\bar B$.
Hence $B=\langle d\rangle$ and we have
$$G=A\times(B\rtimes T)=\langle a\rangle\times(\langle d\rangle\rtimes\langle t\rangle)$$
in which either $|B|\geq3$ in case (b), or $\omega(|B|)<i$ in case (c).
So $G$ is $F$-special by definition and the proof is now complete.
\end{proof}

The more subtle arguments on $F$-special subgroups occur in the following lemmas.

\begin{lem}\label{a1}
The Sylow $2$-subgroup $T$ of $G$ is of order $2$.
\end{lem}
\begin{proof}
We fix some $i\in\{1,\ldots,n\}$ and consider the simple $FN$-module $W_i$
defined in Lemma~\ref{6}.
Take a finite field extension $E/F$ such that $E$ is a splitting field for all subgroups of $G$.
According to~\cite[Theorem 9.21]{I}, we may write the extended $EN$-module
$W_i^E=W_i\otimes_FE$ as
$$W_i^E\cong X_i^1\oplus\cdots\oplus X_i^{m_i},$$
where $X_i^1,\ldots,X_i^{m_i}$ are pairwise non-isomorphic simple $EN$-modules,
and thus of dimension one (as $N$ is abelian), that constitute a single $\Gal(E/F)$-orbit.
Since $W_i=(U_i)_N$ is $M$-invariant,
$M$ permutes this Galois orbit.
Let $S_i$ be the stabilizer of $X_i^1$ in $M$. Then $N\subseteq S_i\subseteq M$.
Clearly $N$ acts on $X_i^1$ by scalar multiplication,
which implies that the commutator subgroup $[N,S_i]$ acts on $X_i^1$ trivially.
From this we may easily deduce that $[N,S_i]$ also centralizes $W_i$,
and by Lemma~\ref{6} again, we see that
$$[N,S_i]\subseteq C_N(W_i)=C_N(V_i)\subseteq C_G(V_i).$$

We claim that there exists some $j\in\{1,\ldots,n\}$ such that $S_j=N$.
In fact, if $S_i>N$ for all $i\in\{1,\ldots,n\}$,
then $S_i\cap\langle t^2\rangle>1$, where $T=\langle t\rangle$ so that $M=N\rtimes\langle t^2\rangle$.
In this case $S_i$ contains the unique involution, say $t'$, in $T$, and we have
$$[N,t']\subseteq\bigcap_{i=1}^n [N,S_i]\subseteq\bigcap_{i=1}^n C_G(V_i)=C_G(V)=1.$$
Hence $t'\in C_T(N)=O_2(G)=1$ by Lemma~\ref{faith}.
This contradiction proves that $S_j=N$ for some $j\in\{1,\ldots,n\}$, as claimed.

Note that $S_j=N$ is also the stabilizer of $X_j^{k}$ in $M$ for each $k=1,\ldots,m_j$.
By Clifford's theorem, all $(X_j^k)^M$ are simple $EM$-modules and hence
$(W_j^M)^E\cong (W_j^E)^M$ is a semisimple $EM$-module.
By~\cite[Problem 9.9]{I}, we see that $W_j^M$ must be a semisimple $FM$-module.
Now the Nakayama Reciprocity (see \cite[Theorem VII.4.5]{H2}) implies that
$$\Hom_{FM}(U_j,W_j^M)\cong\Hom_{FN}((U_j)_N,W_j)\cong\Hom_{FN}(W_j,W_j),$$
which is a field by the Schur's lemma. Thus $U_j$ is of multiplicity $1$
in a composition series of the semisimple $FM$-module $W_j^M$.

Assume that $|T|>2$. Then $|M:N|>1$ and $\dim W_j^M>\dim W_j=\dim U_j$.
In this case we may write $W_j^M\cong U_j\oplus U_j'\oplus\cdots$,
where $U_j'$ is a simple $F M$-module with $U_j'\ncong U_j$.
But then, by the Nakayama Reciprocity again, we have
$$\Hom_{F N}(W_j,(U'_j)_N)\cong\Hom_{F M}(W_j^M,U_j')\neq0.$$
Hence $U_j'$ also lies over $W_j$, which forces $U_j'\cong U_j$ as $F M$-modules by Green's theorem.
This contradiction proves that $|T|=2$, as wanted.
\end{proof}

\medskip
According to Lemma~\ref{a1},
we conclude that $G$ is a nonabelian hyper-dihedral group with $O_2(G)=1$.
It follows that each $F$-special subgroup of $G$
must be of the form $C_m$ or $C_m\times D_{2n}$, where $m,n$ are odd integers
satisfying some number-theoretic conditions.
For the sake of brevity, from now on we set
$$T=\langle t\rangle,\;B=[N,t]\rtimes\langle t\rangle,\;\text{and}\;C=C_N(t).$$
Clearly $N=[N,t]\times C$ and $G=B\times C$ by the Fitting's lemma.
Moreover, we know that $B$ is a dihedral group with $4\nmid |B|$
because $t$ inverts $[N,t]$ and $G$ is not abelian by Lemma~\ref{faith}.
Of course $B$ is also an $F$-special group, which forces $C=Z(G)\neq1$.

To get the final contradiction, we now turn to the quotient groups
$G/C_G(V_i)$ for all $i=1,\ldots,n$.

\begin{lem}\label{a2}
For each $i=1,\ldots,n$, $G/C_G(V_i)$ is always a non-cyclic $F$-special group
and $\omega(V_i)=\omega(G/C_G(V_i))+1$.
\end{lem}
\begin{proof}
Note that $(V_i)_N\cong W_i\oplus W^*_i$ is not a simple $FN$-module
by Lemma~\ref{6} and that $|T|=2$ by Lemma~\ref{a1}.
So $t\not\in C_G(V_i)$ and thus
$$C_G(V_i)=C_N(V_i)=C_N(W_i),$$
which implies that $G/C_G(V_i)$ cannot be cyclic and that
each $W_i$ is a faithful simple $N/C_G(V_i)$-module over $F$.
Now Lemma~\ref{omega} yields that
$$\omega(W_i)=\omega(N/C_G(V_i))=\omega(G/C_G(V_i)).$$
Since $\dim V_i=2\dim W_i$, we have $\omega(V_i)=\omega(W_i)+1$ by definition
and the result follows.
\end{proof}

\begin{lem}\label{a3}
Put $\ell+1=\max\{\omega(V_1),\ldots,\omega(V_n)\}$. Then $\omega(G)=\ell\geq0$.
\end{lem}
\begin{proof}
Fix some $i\in\{1,\ldots,n\}$. By Lemma~\ref{a2}, we see that $\omega(G/C_G(V_i))\leq \ell$
and the equality can hold.
If $p$ is an odd prime dividing $|G|$ such that $\omega(p)>\ell$,
then $p$ does not divide $|G:C_G(V_i)|$.
In this case we may choose an element $x\in G$ of order $p$.
Then $x\in C_G(V_i)$ and hence $x\in\bigcap_{i=1}^n C_G(V_i)=C_G(V)=1$ by Lemma~\ref{faith}.
This contradiction shows that $\omega(G)=\ell$, as desired.
\end{proof}

\begin{lem}\label{a4}
If $\omega(V_i)=\ell+1$ for some $i\in\{1,\ldots,n\}$,
then $G/C_G(V_i)$ is not a dihedral group.
\end{lem}
\begin{proof}
Let $\Lambda$ be the set of those $i\in\{1,\ldots,n\}$
such that $\omega(V_i)=\ell+1$ and that $G/C_G(V_i)$ is dihedral,
and suppose that $\Lambda$ is not empty.
Since $G=B\times C$ and $C=Z(G)$, it follows that $C\subseteq C_G(V_i)$ for each $i\in\Lambda$.
Therefore, the restriction $(V_i)_B$, for all such $i$,
are pairwise non-isomorphic self-dual simple $FB$-modules.
From the fact that $C\neq1$ and $C_G(V)=1$,
we conclude that $\Lambda$ is a proper subset of $\{1,\ldots,n\}$.
So we may fix some $j\in\{1,\ldots,n\}-\Lambda$. Let $X=V_j$.
Then Clifford's theorem implies that $X_B\cong eY$, where $e\geq1$
and $Y$ is a simple $F B$-submodule of $X$.

It suffices to prove that $\omega(Y)<\ell+1$.
Since then $\omega(Y)<\omega(V_i)$ for all $i\in\Lambda$,
it follows that $Y\ncong (V_i)_B$ as simple $FB$-modules.
From this we may easily deduce that $V_B$ cannot be an even-multiplicity $FB$-module
and the minimality of $|G|$ forces $G$ to be the dihedral group $B$.
Hence $G$ itself is $F$-special.
This contradiction tells us that the set $\Lambda$ must be empty,
thus proving the lemma.

To prove the inequality $\omega(Y)<\ell+1$, we need to distinguish two cases
according to the value of $\omega(X)$ by Lemma~\ref{a3}.
Assume first that $\omega(X)<\ell+1$.
Then, since $\dim Y$ divides $\dim X$, it follows that $\omega(Y)\leq\omega(X)<\ell+1$,
and we are done in this case.
Now we assume that $\omega(X)=\ell+1$.
Then, by Lemma~\ref{a2} and the choice of $X$,
we see that $G/C_G(X)$ is $F$-special, but neither cyclic nor dihedral.
So, by definition of $F$-special groups, we may write
$$G/C_G(X)\cong C_a\times D_{2b}$$
where $a>1$ is a $\pi_\ell(q)$-number, $\ell\geq1$, and $b$ is odd with $\omega(b)<\ell$.
From this we may easily conclude that the image in $G/C_G(X)$
of the dihedral group $B$ must be isomorphic with $D_{2b}$.
Let $\bar B=B/C_B(Y)$. Then $C_B(X)=C_B(Y)$ and we obtain
$$\bar B=B/C_B(X)\cong BC_G(X)/C_G(X)\cong D_{2b}.$$
Hence $\omega(\bar B)=\omega(b)<\ell$.
Furthermore, since $Y$ is a faithful simple $F\bar B$-module,
it follows from Lemma~\ref{omega} that $\omega(Y)\leq \omega(\bar B)+1<\ell+1$, as wanted.
\end{proof}

Note that $G=B\times C$, where $B$ is dihedral and $C$ is cyclic.
We need to evaluate $\omega(B)$ and $\omega(C)$.

\begin{lem}\label{a5}
$\omega(B)<\ell$. In particular, we have $\ell\geq1$.
\end{lem}
\begin{proof}
Clearly $\omega(B)\leq\omega(G)\leq\ell$ by Lemma~\ref{a3}.
Assume that $\omega(B)=\ell$. Then we may choose an odd prime divisor $p$ of $|B|$
with $\omega(p)=\ell$ and an element $x\in B$ of order $p$.
Fix some $i\in\{1,\ldots,n\}$. If $\omega(V_i)=\ell+1$,
then, as before, by Lemmas~\ref{a2} and~\ref{a4} we may let $G/C_G(V_i)\cong C_a\times D_{2b}$,
where $a>1$ is a $\pi_\ell(q)$-number, $\ell\geq1$ and $b$ is odd with $\omega(b)<\ell$.
In this case, since $D_{2b}$ is clearly a homomorphism image of $B$ and $p\nmid b$,
it follows that $x\in C_G(V_i)$. In the remaining case where $\omega(V_i)<\ell+1$,
we see that $\omega(G/C_G(V_i))<\ell$ by Lemma~\ref{a2}.
But then $p$ cannot divide $|G/C_G(V_i)|$, which implies that $x\in C_G(V_i)$.
Therefore we have
$$x\in \bigcap_{i=1}^nC_G(V_i)=C_G(V)=1,$$
contradicting Lemma~\ref{faith}.
So $\omega(B)<\ell$, and the proof is complete.
\end{proof}

\begin{lem}\label{a6}
$C$ is a $\pi_\ell(q)$-group.
\end{lem}
\begin{proof}
Suppose that $C$ is not a $\pi_\ell(q)$-group.
Then we may let $C=H\times J$, where $H$ is a $\pi_\ell(q)$-subgroup
and $J$ is a nontrivial $\pi_\ell(q)'$-subgroup of $C$.
By Lemma~\ref{a3} we know that $\omega(G)=\ell$ and hence $\omega(J)<\ell$.
Set $R=B\times H$. Then $G=R\times J$.

Now let $\Lambda=\{i|\,\omega(V_i)=\ell+1,\,1\leq i\leq n\}$.
For each $i\in\Lambda$, by Lemmas~\ref{a2} and ~\ref{a4},
we see that $G/C_G(V_i)$, as an $F$-special group,
is neither cyclic nor dihedral. Hence we may write
$$G/C_G(V_i)\cong C_{a_i}\times D_{2b_i}$$
where $a_i>1$ is a $\pi_\ell(q)$-number, $\ell\geq1$, and $b_i$ is odd with $\omega(b_i)<\ell$.
It is easy to see that each $C_{a_i}$ is a homomorphism image of $C$,
which forces $J\subseteq C_G(V_i)$.
From this we may deduce that $(V_i)_R$, for all $i\in\Lambda$,
are pairwise non-isomorphic self-dual simple $FR$-modules.

On the other hand, since $J\neq1$ and $C_G(V)=1$ by Lemma~\ref{faith},
it follows that $\Lambda$ is a proper subset of $\{1,\ldots,n\}$.
For each $V_j$ with $j\not\in\Lambda$, we have $\omega(V_j)<\ell+1$ by definition.
Let $X_j$ be a simple $F R$-submodule of $V_j$.
Then $(V_j)_R\cong e_jX_j$ for some $e_j\geq1$ by Clifford's theorem.
So $\dim X_j$ divides $\dim V_j$ and $\omega(X_j)\leq\omega(V_j)<\ell+1$.
This implies that $X_j\ncong (V_i)_R$ as simple $F R$-modules
for all $j\not\in\Lambda$ and $i\in\Lambda$.
Thus $V_R$ cannot be an even-multiplicity $F R$-module,
which forces $G=R$ and $J=1$. This contradiction proves the lemma.
\end{proof}

\begin{lem}
$G$ is an $F$-special group, which is the final contradiction.
\end{lem}
\begin{proof}
Recall that $C=C_N(t)\neq1$ and $B=[N,t]\rtimes\langle t\rangle$.
Let $|C|=m$ and $|[N,t]|=n$. By definition, we have $\omega(B)=\omega([N,t])=\omega(n)$.
Now Lemmas~\ref{a5} and~\ref{a6} imply that $m>1$ is a $\pi_\ell(q)$-number with $\ell\geq1$,
and $n>1$ is odd with $\omega(n)<\ell$. Since $N=C\times[N,t]$, it follows that
$G=C\times([N,t]\rtimes\langle t\rangle)\cong C_m\times D_{2n}$.
Therefore $G$ itself is an $F$-special group by definition, which cannot be the case,
thus proving the theorem.
\end{proof}

We are now ready to prove Theorem B from the Introduction.

\medskip
\noindent\emph{Proof of Theorem B.}
The implication ``(1)$\Longrightarrow$(2)" is clear,
and ``(2)$\Longrightarrow$(1)" follows by Theorem~\ref{1}.
What remains is to prove ``(1)$\Longleftrightarrow$(3)".

To do this, we see by definition that an odd-order element $g$ of $G$ is $F$-special
if and only if it lies in some $F$-special subgroup $D$ of $G$
such that the cyclic group $\langle g\rangle$ generated by $g$ is the unique normal $2$-complement of $D$.
According to Theorem~\ref{1}, we may therefore assume further that $G$ itself is an $F$-special group.
So we may let $G=N\rtimes T$ in which $N=\langle g\rangle$ is cyclic of odd order
and $T$ is a $2$-group. Of course $T$ is trivial if $G$ is cyclic.
It suffices to prove that $V$ is an even-multiplicity $FG$-module
if and only if the characteristic polynomial of $g$ on $V$ is a square.

In fact, since each simple $FG$-module is self-dual by Corollary~\ref{special-dual},
it follows that $V$ is an even-multiplicity $FG$-module
if and only if each composition factor of $V$ occurs with even multiplicity,
and also if and only if $V_N\cong 2\Delta$ for some semisimple $FN$-module $\Delta$ by Green's theorem.
Now we assume that
$$V_N\cong k_1W_1\oplus\cdots\oplus k_nW_n,$$
where $W_1,\ldots,W_n$ are pairwise non-isomorphic simple $FN$-modules and $k_i\geq1$ for all $i=1,\ldots,n$.
We write $c_i(x)$ to denote the characteristic polynomial of $g$ on $W_i$
and we may easily deduce that $c_1(x),\ldots,c_n(x)$
are distinct irreducible polynomials over $F$, see~\cite[Lemma 2.1]{J2011}.
Let $c(x)=c_1(x)^{k_1}\cdots c_n(x)^{k_n}$. Then $c(x)$ is the characteristic polynomial of $g$ on $V$.
From this we may conclude that $c(x)$ is a square precisely when all $k_i$ are even,
which is clearly equivalent to the condition that $V_N\cong 2\Delta$ for some $FN$-module $\Delta$.
The proof of Theorem B is now complete.
\qed

\bigskip
Finally, to prove Theorem A we need a lemma essentially due to Dade.

\begin{lem}\label{em}
Let $G$ be a group and let $F$ be a finite field of characteristic $2$.
If $V$ is a symmetric $FG$-module,
then $V$ is hyperbolic if and only if it is an even-multiplicity module.
\end{lem}
\begin{proof}
By the same argument used in the proof of~\cite[Corollary 3.8]{J2011},
we may conclude that the Witt kernel of $V$ is a direct sum of
pairwise non-isomorphic self-dual simple $FG$-modules,
and the result follows by Corollary 2.10 in~\cite{D}.
\end{proof}

\medskip
\noindent\emph{Proof of Theorem A.}
Clear from Lemma~\ref{em} and Theorem B !
\qed

%%%%%%%%%%%%%%%%%%%%%%%%%%%%%%%%%%%%%%%%%%%%%%%%%%%%%%%%%%%%%%%%%%%%%%%%%%%%%%%%%%%%%%%%%%%%%%
\section{Examples}
We will present two examples in this section.
The first one is to show that Theorem A can fail if the $F$-special subgroups $D$ involved
only run over all cyclic subgroups and all extended dihedral subgroups.
This example also tells us that
the set of $F$-special subgroups is ``minimal" in the sense that
one cannot further reduce the hyperbolicity of
symplectic $FG$-modules to some smaller subgroups in general.

\begin{thm}
Let $F$ be a finite filed of characteristic $2$ with $q\geq4$ elements,
let $p,r$ be distinct prime numbers such that $p|(q+1)$ and $r|(q-1)$,
and let $G=C_p\times D_{2r}$. Then the following assertions hold.

{\rm (1)} $G$ is an $F$-special group.

{\rm (2)} There exists a symplectic $FG$-module $V$
such that $V_H$ is a hyperbolic $FH$-module for each proper subgroup $H$ of $G$
but $V$ itself is not hyperbolic.
\end{thm}
\begin{proof}
Note that $p\in\pi_1(q)$ and $r\in\pi_0(q)$, and part (1) follows by definition.

To prove part (2), we may choose a faithful self-dual simple $FG$-module $V$
as guaranteed by Theorem~\ref{dih-dual}.
Then the faithfulness implies that $V$ is not the trivial $FG$-module $F$,
and by a result of Fong~\cite[Theorem VII.8.13]{H2},
we see that $V$ admits a non-degenerate $G$-invariant symplectic form.
Thus $V$ is a symplectic $FG$-module.
Of course, the simple $FG$-module $V$ cannot be hyperbolic.
To complete the proof, it suffices to verify
that $V_H$ is hyperbolic for every maximal subgroup $H$ of $G$.

Let $C=C_p$ and let $B=D_{2r}=N\rtimes\langle t\rangle$,
where $N$ is cyclic of order $r$ and $t$ is an involution.
Since the odd prime $r$ divides $q-1$, we see that
$F$ is a splitting field for the dihedral subgroup $B$.
Also, since $G=C\times B$, we may deduce from~\cite[Corollary 8.13, page 661]{K} that
$$V\cong U\otimes_FW$$
as $FG$-modules, where $U$ is a faithful simple $FC$-module and
$W$ is a faithful simple $FB$-module.
Clearly $C$ and $B$ are $F$-special groups by definition,
and Corollary~\ref{special-dual} implies that both $U$ and $W$ are self-dual modules.
Since groups of odd order have no nontrivial self-dual absolutely irreducible modules
(see~\cite[Remark VII.8.22]{H2}),
it follows that $W_N\cong X\oplus X^*$,
where $X$ is a non-self-dual simple $FN$-module of dimension $1$.
Hence $\dim W=2$, and by~\cite[Problem 9.20]{I} we also have $\dim U=2$.

Note that all of the maximal subgroups of $G$ are $B$, $C\times N$ and $C\times\langle t\rangle$.
From the preceding paragraph we may conclude that $V_B\cong 2W$, $V_C\cong 2U$ and
$$V_{C\times N}\cong U\otimes_F(X\oplus X^*)
=(U\otimes_F X)\oplus(U\otimes_F X)^*.$$
This proves that $V_B$, $V_C$ and $V_{C\times N}$ all are even-multiplicity modules
and hence hyperbolic.
Clearly $V_{C\times\langle t\rangle}$ is hyperbolic if and only if $V_C$
is hyperbolic, and the result follows.
\end{proof}

The second example is to demonstrate that Theorem A can also fail
if all the extended dihedral subgroups involved
are replaced by the usual dihedral subgroups,
even in the case where $V$ is a faithful $FG$-module and
$F$ is a splitting field for all subgroups of $G$.

\begin{thm}
Let $F$ be a finite field with $q$ elements,
where $q\geq4$ is a power of $2$.
Let $p$ be any prime divisor of $q-1$ and let $G=\tilde D_{4p}$.
Then $F$ is a splitting field for all subgroups of $G$, and there exists a faithful symplectic
$FG$-module $V$ such that $V_H$ is hyperbolic for each proper subgroup $H$ of $G$
but $V$ itself is not hyperbolic.
\end{thm}
\begin{proof}
That $F$ is a splitting field for all subgroups of $G$ is clear.
Since $G/O_2(G)$ is the usual dihedral group $D_{2p}$,
by Theorem~\ref{dih-dual} we may take $U$ to be a faithful self-dual simple
$F[G/O_2(G)]$-module and lift it to an $FG$-module in an obvious way.
Then the Fong's theorem cited above implies that $U$ is a symplectic $FG$-module.
In particular, the $F$-dimension of $U$ must be even.
To get a faithful $FG$-module, we consider the regular module $FG$.
By~\cite[Exercise 24, page 120]{H2},
we may choose a $G$-invariant non-degenerate symplectic form on $FG$
so that it becomes into a faithful symplectic $FG$-module.
It is easy to see that any nontrivial element in $G$ of odd order
has the characteristic polynomial on $FG$ of the form $(x^p-1)^4$,
and hence Theorem A implies that the regular module $FG$ is hyperbolic.

Now let $V=FG\bot U$ be the orthogonal direct sum of these two symplectic $FG$-modules.
Then $V$ clearly cannot be hyperbolic as a symplectic $FG$-module.
However, for every proper subgroup $H$ of $G$,
since $U$ is of even dimension,
the characteristic polynomial on $V$ defined by every element of $H$ of odd order
must be a square. By Theorem A again, we see that
$U_H$ is hyperbolic and the proof is complete.
\end{proof}

%%%%%%%%%%%%%%%%%%%%%%%%%%%%%%%%%%%%%%%%%%%%%%%%%%%%%%%%%%%%%%%%%%%%%%%%%%%%%%%%%%%%%%
\section{Applications}
In the final section we will present some immediate applications
of our main theorems to the real-valued Brauer characters at prime $2$,
the automorphisms of self-dual binary codes,
and the Witt groups of symplectic modules.
In the sequel we will demonstrate how Theorem A can be used
to reduce certain problems on monomial characters and
on the Glauberman-Isaacs character correspondences
to the same problems for $F$-special subgroups.

\subsection{Real-valued irreducible Brauer characters at prime $2$}
In order to study the Frobenius-Schur indicator of a real-valued
irreducible complex character of a finite group,
R. Gow obtained the following result~\cite[Theorem 1.1]{G},
which is a powerful tool.

\begin{thm}[Gow]
Let $\chi$ be a real-valued irreducible character of a finite group $G$.
Then there is a subgroup $H$ of $G$ and a real-valued irreducible character $\theta$ of $H$
for which $[\theta,\chi_H]$ is odd.
$H$ can be taken either to be a Sylow $2$-subgroup of $G$
or to have the form $AU$, where $A=\langle h\rangle$
is a cyclic subgroup of odd order generated by a real non-identity element $h$
and $U$ is a Sylow $2$-subgroup of $C^*(h)$ defined by
$$C^*(h)=\{x\in G\,|\,h^x=h\;\text{or}\;h^{-1}\}.$$
\end{thm}

Now we give a similar version for $2$-Brauer characters,
which we think is also useful in the study of real-valued irreducible Brauer characters.

\begin{thm}\label{real-valued Brauer}
Let $G$ be a group and let $\varphi\in\IBr_2(G)$
be a nontrivial real-valued irreducible Brauer character of $G$ at prime $2$.
Then there exist a subgroup $H$ and a real-valued irreducible Brauer character $\theta\in\IBr_2(H)$
for which the multiplicity of $\theta$ in $\varphi_H$ is odd.
Moreover, $H$ can be chosen as a semi-direct product $\langle h\rangle\rtimes\langle t\rangle$,
where $h$ is a real element of odd order and $t$ is a $2$-element
inverting $h$. In particular, the inner product $[\varphi_{\langle h\rangle},\lambda]$
of complex characters is odd for some $\lambda\in\Irr(\langle h\rangle)$.
\end{thm}
\begin{proof}
Let $F$ be a finite field of characteristic $2$
such that it is a splitting field for all subgroups of $G$.
Then $\varphi$ is the Brauer character afforded by some self-dual simple $FG$-module $V$.
By Theorem B, there exists an extended dihedral group $H\leq G$
such that $V_H$ cannot be an even-multiplicity module,
that is, by definition, we have some self-dual simple $FH$-module $W$
with odd multiplicity in $V_H$.
Let $\theta$ be the Brauer character of $H$ afforded by $W$.
Then $\theta$ is a real-valued irreducible Brauer character of $H$
with odd multiplicity in $\varphi_H$.
Furthermore, by definition of extended dihedral groups,
we may write $H=\langle h\rangle\rtimes\langle t\rangle$,
where $h\neq1$ is a real element of odd order and $t$ is some $2$-element
such that $h^t=h^{-1}$.
Finally, since $h$ is $F$-special, it follows from Theorem B
that $V_{\langle h\rangle}\not\cong 2\Delta$ for any $F\langle h\rangle$-module $\Delta$,
which amounts to saying that $[\varphi_{\langle h\rangle},\lambda]$ is odd
for some $\lambda\in\Irr(\langle h\rangle)$. The proof is now complete.
\end{proof}

\subsection{Automorphisms of self-dual binary codes}
A. G\"unther and G. Nebe in~\cite{GN} studied the automorphism group of self-dual binary codes.
Recall that a \emph{linear binary code} of length $n$ is merely a subspace $C\subseteq\mathbb{F}_2^n$,
and $C$ is said to be \emph{self-dual} if $C$ coincides with its dual $C^\bot$ defined by
$$C^\bot=\{v\in\mathbb{F}_2^n\,|\,v\cdot c=\sum_{i=1}^nv_ic_i=0\;\;\text{for all}\;\;c\in C\}.$$
The automorphism group $\Aut(C)$ of $C$ consist of those permutations $\sigma\in S_n$
(the symmetric group of degree $n$) that leave $C$ invariant by permuting coordinates, i.e.
$$\Aut(C)=\{\sigma\in S_n\,|\,C\sigma=C\}.$$
One of the main results that they proved is the characterization of permutation groups
that act on a self-dual code (see~\cite[Theorem 2.1]{GN}).

\begin{thm}[G\"unther-Nebe]
Let $G\leq S_n$. Then there exists a self-dual code $C\subseteq\mathbb{F}_2^n$
with $G\leq\Aut(C)$ if and only if every self-dual simple $\mathbb{F}_2G$-module
occurs in the $\mathbb{F}_2G$-module $\mathbb{F}_2^n$ with even multiplicity.
\end{thm}

In our words, this theorem amounts to saying that
the symmetric $\mathbb{F}_2G$-module $\mathbb{F}_2^n$
(with respect to the standard scalar product) is hyperbolic if and only if
it is an even-multiplicity module.
As an immediate application of our Theorem A,
we will present a more effective criterion.
In particular, we characterize the full automorphism group of such self-dual codes.

\begin{thm}\label{code}
Let $G\leq S_n$. Then there exists a self-dual code $C\subseteq\mathbb{F}_2^n$
with $G\leq\Aut(C)$ if and only if $G$ satisfies the condition:
For every $\mathbb{F}_2$-special element $\sigma\in G$,
the number of cycles of equal size in the cyclic decomposition of $\sigma$ is always even.

In particular, $G=\Aut(C)$ for some self-dual code $C\subseteq\mathbb{F}_2^n$
if and only if $G$ is maximal among subgroups of $S_n$ satisfying this condition.
\end{thm}
\begin{proof}
Let $e_1,\ldots,e_n$ be the standard basis of $\mathbb{F}_2^n$.
Then $S_n$ permutes the set $\{e_1,\ldots,e_n\}$ in an obvious way.
For any $\sigma\in S_n$, let $O_1,\ldots,O_r$ be the orbits of
the cyclic group $\langle\sigma\rangle$ on this set
with $|O_i|=k_i,1\leq i\leq r$.
Then we may deduce that the characteristic polynomial of $\sigma$
on the subspace $\mathbb{F}_2[O_i]$ spanned by orbit $O_i$ is $x^{k_i}-1$.
It follows that the characteristic polynomial of $\sigma$ on the whole space $\mathbb{F}_2^n$
equals to $\prod_{i=1}^r(x^{k_i}-1)$, which is a square if and only if each $k_i$ occurs even times.
Clearly $k_1,\ldots,k_r$ are precisely the sizes of cycles of $\sigma$,
and the result follows by Theorem A.
\end{proof}

It is easy to see that Theorem~\ref{code} also holds if we replace
the ground field $\mathbb{F}_2$ by any finite field $F$ of characteristic two.

\subsection{The Witt group of symplectic modules}
Fixing a finite field $F$ and a finite group $G$, we consider all the symplectic $FG$-modules.
Recall that if $(V,f)$ and $(V',f')$ are symplectic $FG$-modules, the orthogonal sum is defined by
$$(V,f)\,\bot\,(V',f')=(V\oplus V',\,f\bot f'),$$
where $(f\bot f')(v+v',w+w')=f(v,w)+f'(v',w')$ for $v,v'\in V,w,w'\in W$.
As usual, we say that $(V,f)$ and $(V',f')$ are \emph{Witt equivalent}
if $(V,f)\bot(V',-f')$ is hyperbolic.
Denote by $[(V,f)]$ the Witt equivalence class of $(V,f)$
and let $\mathcal{W}_s(FG)$ be the Witt group consist of
the Witt equivalence classes of symplectic $FG$-modules with multiplication
$$[(V,f)]\,\bot\,[(V',f')]=[(V,f)\,\bot\,(V',f')],$$
which is well-defined.
A simple fact is that two symplectic $FG$-modules $V$ and $V'$ are Witt equivalent
if and only if they have the same Witt kernel (up to $G$-isometry).
Of course, all of these discussions are quite classical, see~\cite{Sc}, for example.

As an application of Theorem A, we prove the following.

\begin{thm}\label{witt}
Let $F$ be a finite field of characteristic $2$ and let $G$ be a group.
For any two symplectic $FG$-modules $V$ and $V'$, the following conditions are equivalent:

{\rm(1)} $V$ and $V'$ are Witt equivalent, that is, $[V]=[V']\in\mathcal{W}_s(FG)$.

{\rm(2)} $V_D$ and $V'_D$ are Witt equivalent for every $F$-special subgroup $D$ of $G$.

{\rm(3)} For each $F$-special element $g$ of $G$,
the product of the characteristic polynomials of $g$ on $V$ and on $V'$ is a square.
\end{thm}
\begin{proof}
Let $D_1,\ldots, D_n$ be all the $F$-special subgroups of $G$.
Then it is easy to verify that the restriction
$${\rm Res}:\mathcal{W}_s(FG)\to \bigoplus_{i=1}^n\mathcal{W}_s(FD_i),\;\;
[V]\mapsto [V_{D_1}]+\cdots+[V_{D_n}]$$
is well-defined and gives rise to a group homomorphism.
By Theorem A, we see that this map is injective,
which proves the equivalence  ``(1)$\Longleftrightarrow$(2)".
Also, since $V$ and $V'$ are Witt equivalent if and only if
$V\bot V'$ is hyperbolic, it follows that
the equivalence ``(1)$\Longleftrightarrow$(3)" by Theorem A again.
\end{proof}

Clearly, Theorem~\ref{witt} is also valid for the Witt group of symmetric $FG$-modules
if we replace the word ``symplectic" with ``symmetric".

%%%%%%%%%%%%%%%%%%%%%%%%%%%%%%%%%%%%%%%%%%%%%%%%%%%%%%%%%%%%%%%%%%%%%%%%%%%%%%%%%%%%%%

%%%-------------------------------------------------------------------
\end{document}